\documentclass[11pt,twoside]{amsart}

\usepackage{amsfonts}

\title{Degenerations for derived categories}
\author{Bernt Tore Jensen}
\address{Bernt Tore Jensen,\newline Department of Pure Mathematics,\newline
University of Leeds,\newline
Leeds LS2 9JT,\newline UK}
\email{bjensen@maths.leeds.ac.uk}
\author{Xiuping Su}
\address{Xiuping Su, Alexander Zimmermann\newline
Universit\'e de Picardie,\newline Facult\'e de Math\'ematiques et
LAMFA (UMR 6140 du CNRS), \newline 33 rue St Leu,\newline F-80039
Amiens Cedex 1,\newline France} \email{xiuping.su@u-picardie.fr,
alexander.zimmermann@u-picardie.fr}
\author{Alexander Zimmermann}
\date{December 4, 2003; revised September 24, 2004}
\dedicatory{To Claus Michael Ringel for his sixtieth birthday}

\newtheorem{Theo1}{{Theorem}}
\newtheorem*{Theo2}{{Theorem}}
\newtheorem{Lemma1}{{Lemma}}
\newtheorem{Def1}{{Definition}}
\newtheorem{Prop1}[Lemma1]{{Proposition}}
\newtheorem{Claim1}[Lemma1]{{Claim}}
\newtheorem{Rem1}{{Remark}}[section]
\newtheorem{Cor1}[Lemma1]{{Corollary}}
\newtheorem{Ex1}[Rem1]{{Example}}

\newenvironment{Lemma}{\begin{Lemma1}}{\end{Lemma1}}

\newenvironment{Prop}{\begin{Prop1}}{\end{Prop1}}
\newenvironment{Rem}{\begin{Rem1}\rm}{\end{Rem1}}
\newenvironment{Theorem}{\begin{Theo1}}{\end{Theo1}}
\newenvironment{Theorem2}{\begin{Theo2}}{\end{Theo2}}
\newenvironment{Remark}{\begin{Rem1}\em}{\end{Rem1}}
\newenvironment{Cor}{\begin{Cor1}}{\end{Cor1}}

\newenvironment{Example}{\begin{Ex1}\em}{\end{Ex1}}

\newcommand{\dar}{\downarrow}
\newcommand{\lra}{\longrightarrow}

\newcommand{\ra}{\rightarrow}
\newcommand{\sdp}{\times\kern-.2em\vrule height1.1ex depth-.05ex}
\newcommand{\epi}{\lra \kern-.8em\ra}

\newcommand{\N}{{\mathbb N}}

\newcommand{\Z}{{\mathbb Z}}

\newcommand{\dickebox}{{\vrule height5pt width5pt depth0pt}}
\newcommand{\Scomp}{comproj}

 \normalbaselines
\begin{document}

\begin{abstract}
We propose a theory of degenerations for derived module categories,
analogous to degenerations in module varieties for module categories. In
particular we define two types of degenerations, one algebraic
and the other geometric. We show that these are equivalent,
analogously to the Riedtmann-Zwara theorem for module varieties.
Applications to tilting complexes are given, in particular that any
two term tilting complex is determined by its graded module structure.
\end{abstract}

\maketitle

\section*{Introduction}

Geometrical methods were introduced in representation theory
of finite dimensional algebras in order to parameterize
 possible module structures on a given vector space by algebraic
varieties. These varieties carry an action of a reductive
algebraic group $G$ such that the orbits correspond to isomorphism
classes of modules. One says that a module $M$ degenerates to $N$
if $N$ is in the closure $\overline{G\cdot M}$ of the orbit of $M$
under the $G$-action, and in this case one writes $M\leq N$.
Riedtmann defined in \cite{Riedtmann} a relation $\leq_{alg}$ by
setting $M\leq_{alg} N$ if there is a module $Z$ and a short exact
sequence $$0\lra N \lra M\oplus Z \lra Z \lra 0$$ of $A$-modules.
She showed that $M\leq_{alg}N$ implies $M\leq N$ and in
\cite{Zwaramain} Zwara proved that $M\leq N$ implies
$M\leq_{alg}N$.

Since the derived category became a powerful tool in representation theory
it seems desirable to study derived categories from such a
geometric point of view.
De Concini and Strickland \cite{deConciniStrickland} studied
geometric properties of
varieties of bounded complexes of free modules. For a
finite dimensional algebra $A$
Huisgen-Zimmermann and Saorin \cite{SaorinHuisgenZimmermann}
defined an affine variety which parameterizes bounded complexes of
$A$-modules. For this variety no group action seems available so that
the quasi-isomorphism classes correspond to orbits under the
action.
Bekkert and Drozd studied in \cite{DrozdBekkert}
minimal right bounded complexes of projective modules,
where quasi-isomorphism
is the same as homotopy equivalence. There homotopy equivalence classes
are obtained as orbits of an action of a group, however Bekkert
and Drozd did not study the topology of their space and in particular
they did not study degeneration.

The purpose of the present paper is to define and to study a
geometric structure on a set of right bounded complexes of projective
modules and to show a result analogous to the result of Zwara and
Riedtmann.
More precisely, we define a topological space $\Scomp^{\underline d}$
parameterizing right bounded complexes of projective modules
depending on a dimension array ${\underline d}$
replacing the dimension vector for module varieties.
This topological space is a projective limit of affine varieties and a
projective limit $G$ of affine algebraic groups is acting on it.
The $G$-orbits correspond to quasi-isomorphism classes of right bounded
complexes of projective modules.
For two right bounded complexes $M$ and $N$ we define $M\leq_{\Delta}N$
if there is a complex $Z$ and a distinguished triangle
$$N\lra M\oplus Z\lra Z\lra N[1]\;.$$
For two right bounded complexes $M$ and $N$ in $\Scomp^{\underline d}$
we say $M\leq_{top}N$ if $N\in\overline{G\cdot M}$.
Our main result is the following.

\begin{Theorem2}
Let $A$ be a finite dimensional
$k$-algebra over an algebraically closed field $k$ and let
$N$ and $M$ be complexes in the bounded derived category of $A$-modules
$D^b(A)$.
Then, there is a dimension array $\underline d$ so that
$N$ and $M$ belong to $\Scomp^{\underline d}$
and moreover $M\leq_{\Delta} N$ if and only if $M\leq_{top} N$.
\end{Theorem2}

Using $\leq_{alg}$ and $\leq_{\Delta}$ we show that for two $A$-modules
$M$ and $N$ one can choose a dimension array
$\underline d$ so that the module $M$ degenerates to $N$ in the
module variety if and only if the projective resolution of $M$
degenerates to the projective resolution of $N$ in $\Scomp^{\underline d}$.
To illustrate how the topology of $\Scomp^{\underline d}$ can be used
we show that a partial
two-term tilting complex is determined, up to isomorphism, by its
structure as a graded module. We give an example showing that this is
not true for longer tilting complexes.

The paper is organized as follows. In Section~\ref{basics} we
define the variety $\Scomp^{\underline d}$, define a group acting on it,
and show some basic properties. In Section~\ref{algebricimpliestopological}
we define $\leq_{\Delta}$ and show that $\leq_{\Delta}$ implies the
topological degeneration for two complexes with bounded homology.
In Section ~\ref{converse} we show the converse. Section~\ref{consequences}
finally develops consequences for complexes without self-extensions.

\section{General definitions and elementary properties}
\label{basics}

Let $A$ be a finite dimensional algebra over an algebraically
closed field $k$. Let $mod(A,d)$ denote the affine variety of
$d$-dimensional $A$-modules. The general linear group $Gl_d(k)$
acts on $mod(A,d)$ by change of basis and the orbits correspond to
the isomorphism classes of $d$-dimensional modules.

Let $P_1,\dots,P_l$ be a complete set of projective indecomposable
$A$-modules one in each isomorphism class. For an element
$d=(d^1,\dots ,d^l) \in \mathbb{N}^l$ let $\alpha(d)$ be defined by
$\bigoplus_{j=1}^lP_j^{d^j}\in mod(A,\alpha(d))$.

For every sequence $\underline d:\Z\lra \N^{l}$ for which there is
an $i_0\in\Z$ with $d_i=(0,\dots,0)$ for $i\leq i_0$ we define
$comp(A,\underline\alpha(\underline d))$ to be the subset of
$$\left(\prod_{i\in\Z}mod(A,\alpha(d_i))\right)\times
\left(\prod_{i\in\Z}Hom_k(k^{\alpha(d_i)},k^{\alpha(d_{i-1})})\right)$$
consisting of elements $((M_i)_{i\in \Z},(\partial_i)_{i\in \Z})$
with the properties that $\partial _i$ is an $A$-homomorphism when
viewed as a map from $M_i$ to $M_{i-1}$ and $\partial_i \partial_{i-1} = 0$.

The group $\prod_{i\in \Z} Gl_{\alpha(d_i)}$ acts on
$comp(A,\underline\alpha(\underline d))$ by change of basis
and the orbits correspond to isomorphism classes of complexes.

We have a projection $\pi_M:comp(A,\underline\alpha(\underline d))
\lra \prod_{i\in\Z}mod(A,\alpha(d_i))$ and we define
$$\Scomp^{\underline d}:= \pi_M^{-1}(\prod_{i\in\Z}
\bigoplus_{j=1}^lP_j^{d_i^j})\;.$$

We say that $\underline{d}$ is bounded if there is an $i_1\in\Z$ with
$d_i=(0,\dots,0)$ for $i\geq i_1$. In this case we identify
$comp(A,\underline\alpha (\underline d))$ with the affine variety of
bounded complexes defined by Huisgen-Zimmermann and Saorin in
\cite{SaorinHuisgenZimmermann}, in particular it has the
Zariski topology. Also $\Scomp^{\underline d}$ is then an affine
variety, being a closed subset of $comp(A,\underline\alpha(\underline d))$.

Naive truncation on the left induces surjective
morphisms of varieties
\begin{eqnarray*}
\varphi_n:\Scomp^{\underline d_n}&\lra& \Scomp^{\underline d_{n-1}}\\
(\prod_{i\leq n}M_i,\prod_{i\leq n}\partial_i)&\mapsto&
(\prod_{i\leq n-1}M_i,\prod_{i\leq n-1}\partial_i)
\end{eqnarray*}
and similarly surjective maps
\begin{eqnarray*}
\pi_n:\Scomp^{\underline d}&\lra& \Scomp^{\underline d_{n}}\\
(\prod_{i\in\Z}M_i,\prod_{i\in\Z}\partial_i)&\mapsto&
(\prod_{i\leq n}M_i,\prod_{i\leq n}\partial_i)\;.
\end{eqnarray*}


We give $\Scomp^{\underline d}$ the weak topology with
respect to the maps $\{\pi_{n}\}$.
So, the open sets in $\Scomp^{\underline d}$ are of the form
$U=\bigcup_{n\geq n_0}\pi_{n}^{-1}(U_{n})$ for open sets $U_{n}$ in
$\Scomp^{\underline d_{n}}$ and $n_0\in\Z$.
Similarly,
the closed sets in $\Scomp^{\underline d}$ are of the form
$C=\bigcap_{n\geq n_0}\pi_{n}^{-1}(C_{n})$ for closed sets $C_{n}$ in
$\Scomp^{\underline d_{n}}$ and an $n_0\in\Z$.
Note that $\Scomp^{\underline d}$ is the projective limit of the
varieties $\Scomp^{\underline d_n}$ in the category of topological
spaces.

The group $$G:=\prod _{i\in \Z}
Stab_{Gl_{\alpha(d_i)}}(\bigoplus_{j=1}^lP_j^{d_i^j})\cong
\prod _{i\in \Z}Aut _A(\bigoplus_{j=1}^lP_j^{d_i^j})$$
acts on the space $\Scomp^{\underline d}$ by
conjugation and the orbits correspond to isomorphism classes
of complexes of projective $A$-modules. The action of $G$
on $\Scomp^{\underline d}$ induces naturally an action
of $G$ on $\Scomp^{\underline d_n}$
for all $n$ such that $\pi _n$ and $\varphi _n$ are $G$-equivariant maps.

We see that $G$ is a connected algebraic group if $\underline d$ is bounded
since the endomorphism ring is a linear space, hence irreducible, and
the automorphism group is an open dense subvariety.
Moreover, the action of $G$ is the action of a connected algebraic
group on an affine variety if $\underline{d}$ is bounded.

The following lemma is well known to the experts, but we could not
find a reference and include a proof below. We do not require
the field to be algebraically closed for the remainder of this section.

\begin{Lemma}\label{homotopyequalsisomorphism}
Let $X=(\bigoplus_{i\in \Z}Q_i,\bigoplus_{i\in\Z}\partial_i^X)$
and $Y=(\bigoplus_{i\in\Z}Q_i,\bigoplus_{i\in\Z}\partial_i^Y)$
be two right bounded complexes of projective $A$-modules with the
same homogeneous
components $Q_i$ in each degree $i\in\Z$. Then, $X$ is isomorphic to $Y$
if and only if $X$ is homotopy equivalent to $Y$.
\end{Lemma}

Proof.
If $X$ is isomorphic to $Y$ in the category of complexes, then clearly
$X$ is homotopy equivalent to $Y$. So, suppose that $X$ is homotopy
equivalent to $Y$, that is there is a mapping of complexes
$\varphi:X\lra Y$ and a mapping $\psi:Y\lra X$ of complexes
so that there is a map $h$ of degree $1$ so that
$\varphi\psi-id_X=h\partial^X+\partial^Xh$ and likewise
there is an $h'$ with $\psi\varphi-id_Y=h'\partial^Y+\partial^Yh'$.
We shall show that $X\simeq X'\oplus N_X$ where
$im(\partial^X|_{X'})\subseteq rad(X')$ and $N_X$ is contractible,
and likewise for $Y$.

Suppose for the moment that this is shown. Then,
since $N_X$ and $N_Y$ are contractible, we get $X'$ and $Y'$ are
quasi-isomorphic, and therefore, since both are
right bounded complexes of projective modules, homotopy equivalent.
Once we can show that then $X'$ and $Y'$ are isomorphic as complexes,
then also $N_X$ and $N_Y$ are isomorphic. Indeed, $N_X$ and $N_Y$
are isomorphic as graded modules. Now, since $N_X$ and $N_Y$ are
contractible, they are both isomorphic to a direct sum of copies of
shifted copies of complexes of the form
$$\dots\lra 0\lra M\stackrel{\simeq}{\lra} M\lra 0\lra \dots.$$
Comparing the direct factors, and using that $N_X$ and $N_Y$ are isomorphic
as graded modules, one sees that $N_X\simeq N_Y$ as complexes.
So, we suppose for the moment in the statement of the lemma that
$im(\partial^X)\subseteq rad(X)$. But then,
$\varphi\psi-id_X=h\partial^X+\partial^Xh$, and therefore,
$(\varphi\psi-id_X)(X_i)\subseteq rad(X_i)$ for any degree $i$.
Nakayama's Lemma implies that $\varphi\psi$ is invertible. Likewise
$\psi\varphi$ is invertible. Hence, $\varphi$ is an isomorphism.

We need to show that $X\simeq X'\oplus N_X$ for a contractible $N_X$ and
a complex $X'$ with $im(\partial^X|_{X'})\subseteq rad(X')$.
Since $\partial^X(rad(X_i))\subseteq rad(X_{i-1})$, the complex $X$
induces a complex
$(X/rad(X),\overline{\partial}^X)$ of semisimple modules.
Let $m$ be the smallest degree such that $X_{m-1}$ is non-zero. Denote
$\overline{X}:=X/rad(X)$. If
$$0\neq \overline{\partial}^X_m:\overline{X_m}\lra \overline{X}_{m-1},
\mbox{ then }
\overline{X}_{m-1}\simeq \overline{X}_{m-1}'\oplus \overline{X}_{m-1}''
$$
so that $\overline{\partial}^X_m:\overline{X}_m\lra \overline{X}_{m-1}''$
is surjective. But then, $\partial^X_m$ is also surjective onto
the projective cover $X_{m-1}''$ of $\overline{X}_{m-1}''$. Since
$X_{m-1}''$ is projective, there is a splitting
$\sigma_{m-1}: X_{m-1}''\lra X_m$ of $\partial^X_m$ and therefore,
$X_m\simeq X_m'\oplus X_{m-1}''$ and
$X_{m-1}\simeq X_{m-1}'\oplus X_{m-1}''$
so that $\partial^X_m$ is transformed by these isomorphisms into
$\left(\begin{array}{cc}{\partial_m^X}'&0\\
0&id_{X_{m-1}''}\end{array}\right)$.
Now, by induction define $N_X:=X''$ and one gets $X\simeq X'\oplus N_X$
and $\partial^X|_{X'}$ induces the $0$-mapping modulo the radical. This is
tantamount to saying that $im(\partial^X|_{X'})\subseteq rad(X')$.
\phantom{x}\hfill\dickebox

\bigskip

As a consequence of the lemma we see that the orbits in
$\Scomp^{\underline d}$ under the action of G correspond to
homotopy equivalence classes, or equivalently quasi-isomorphism classes,
of right bounded complexes of projective
modules with fixed dimension array ${\underline d}$. Note however
that ${\underline d}$ is not preserved under quasi-isomorphism.

\begin{Lemma}\label{equalizingdimarrays}
Let $M$ and $N$ be right bounded
complexes of finitely generated projective modules.
Then, there is a dimension array $\underline d$ and homotopy equivalent
complexes $M\simeq M'$ and $N\simeq N'$ so that $M',N'\in
\Scomp^{\underline d}$.
\end{Lemma}

Proof.
Let $n$ be the smallest degree such that the homogeneous component of $M$ or $N$ is nonzero.
$$M'_m:=\bigoplus_{j=n}^mM_j\oplus\bigoplus_{j=n}^{m}N_j
\mbox{ and } N'_m:=\bigoplus_{j=n}^mN_j\oplus\bigoplus_{j=n}^{m}M_j$$
where the differential $d_{M'}$ is chosen to be $d_M$ on $M_m$,
and the differential $d_{N'}$ is chosen to be $d_N$ on $N_m$.
Moreover, $$d_{M'}|_{M_k}=
\begin{cases}id&\mbox{if }m-k>0\mbox{ is
even}\\0&\mbox{else}\end{cases}\mbox{
whereas  }d_{M'}|_{N_k}=
\begin{cases}id&\mbox{if }m-k\mbox{ is odd}\\0&\mbox{else}\end{cases}$$
Define the differential on $N'$ likewise, and get this way that
$M'_m\simeq N'_m$ for all $m$, and also $M\simeq M'$ as well as $N\simeq
N'.$
\hfill\dickebox

\bigskip

We define for a complex $X$ the complex $X[1]$ shifted by one
degree to the left by $(X[1])_m:=X_{m-1}$ and
$\partial^{X[1]}_m:=-\partial^X_{m-1}$.

\section{Algebraic relation implies topological relation}
\label{algebricimpliestopological}

Let $A$ be an algebra over an algebraically closed field $k$,
let $D^-(A)$ be the derived category of right bounded
complexes of finitely generated
$A$-modules, and let $D^b(A)$ be its full subcategory
formed by bounded complexes of $A$-modules.
Let $K^-(A)$ be the homotopy category of
right bounded complexes of finitely generated
projective $A$-modules and $K^{-,b}(A)$ the image of
$D^b(A)$ in $K^-(A)$ under the equivalence $K^{-}(A)\simeq D^-(A)$.
Concerning conventions for derived categories we shall
follow \cite{derbuch}.

For any $X$ and $Y$ in $D^-(A)$ let $X\leq_{\Delta} Y$ if there
is a distinguished triangle $$ Y\lra X\oplus Z\lra Z\lra Y[1]$$ for
an object $Z$ in $D^-(A)$.

On the topological side we define a relation $\leq_{top}$ by
$$X\leq_{top}Y:\Leftrightarrow Y\in \overline{G\cdot X}$$
for $X,Y\in \Scomp^{\underline d}$.

We denote by $\underline{dim}(X)$ the dimension array of a
complex $X\in K^{-}(A)$.

Observe that if $X$ and $Y$ are in $D^b(A)$, then $X\leq_{\Delta} Y$
implies $[X]=[Y]$ in $K_0(D^b(A))$.

\begin{Theorem}\label{adherence}
Let $M$ and $N$ be right bounded complexes of finitely generated
projective modules with the same dimension array $\underline{d}$.
Then, $M\leq_{\Delta} N$ implies $M\leq_{top}N$ in
$\Scomp^{\underline{d}}$.
\end{Theorem}

Proof.
Let $U$ be a subset of $\Scomp^{\underline d}$. We show that
${\overline U} = \bigcap _n {\pi_n}^{-1}(\overline{\pi_n(U)})$.
The inclusion $\subseteq$ is obvious. Let $C$ be a closed subset of
$\Scomp^{\underline d}$ with $U\subseteq C$.
Then $C=\bigcap _n {\pi_n}^{-1}(C_n)$
for closed subsets $C_n\subseteq \Scomp^{\underline d_n}$. Hence $U\subseteq
{\pi_n}^{-1}(C_n)$ and so $\overline{\pi_n(U)}\subseteq C_n$
for every $n$, which proves the other inclusion.

Now, if one can prove that whenever $M\leq_{\Delta} N$, then
$\pi_n(M)\leq_{\Delta} \pi_n(N)$, and moreover, if this implies
that $\pi_n(N)\in \overline{G\cdot \pi_n(M)}$, then
by the above, $N\in \overline{G\cdot M}$. This means, that
$M$ degenerates to $N$ in the topological sense.

We still have to show that if $M\leq_{\Delta} N$ then
$\pi_n(N)\in \overline{G\cdot \pi_n(M)}$. We shall use the very same proof
as in the module case by Riedtmann \cite{Riedtmann}.
Let $M\leq_{\Delta} N$. Then, there is a complex $Z$
of projective modules so that
$$N\lra M\oplus Z\lra Z\lra N[1]$$ is a distinguished triangle.
This implies that $$Z[-1]\lra N\lra M\oplus Z\lra (Z[-1])[1]$$
is a distinguished triangle. Hence, $M\oplus Z\simeq cone(Z[-1]\lra N)$
in the homotopy category. Now, we use that the dimension array of
$N$ and of $M$ coincide. Indeed,
$$\underline{dim}(cone(Z[-1]\lra N))=
\underline{dim}(Z)+\underline{dim}(N)=
\underline{dim}(Z)+\underline{dim}(M)=\underline{dim}(M\oplus Z)\;.$$
Hence, $cone(Z[-1]\lra N)\simeq M\oplus Z$ in the
category of complexes and so there is a sequence
$$0\lra N\stackrel{(\phi,\alpha)}{\lra} Z\oplus M
\stackrel{\beta\choose \psi}{\lra} Z\lra 0$$
which is exact in the category of complexes. This shows at once that
$M\leq_{\Delta}N$ implies $\pi_n(M)\leq_{\Delta} \pi_n(N)$ for any $n$.

The first assertion is that $\beta$ is invertible if and only
if $\alpha$ is invertible and in this case,
$$0\lra N\stackrel{(\phi,\alpha)}{\lra} Z\oplus M
\stackrel{\beta\choose \psi}{\lra} Z\lra 0$$
is isomorphic to
$$0\lra N\stackrel{(0,\alpha)}{\lra}Z\oplus M
\stackrel{\beta\choose 0}{\lra}Z\lra 0$$ and therefore $N\simeq M$.
Indeed, we get an isomorphism of exact sequences
$$
\begin{array}{ccccccccc}
0&\lra& N&\stackrel{(\phi,\alpha)}{\lra} &Z\oplus M&
\stackrel{\beta\choose \psi}{\lra}& Z&\lra& 0\\
&&\|&&\phantom{\mbox{\scriptsize$\left(\begin{array}{cc}id_Z&0\\0&id_M
\end{array}\right)$}}
\dar\mbox{\scriptsize$\left(\begin{array}{cc}id_Z&0\\
\psi\beta^{-1}&id_M
\end{array}\right)$}&&\|\\
0&\lra& N&\stackrel{(0,\alpha)}{\lra}&Z\oplus M&
\stackrel{\beta\choose 0}{\lra}&Z&\lra& 0
\end{array}
$$
and likewise for $\alpha$ invertible.


For any $t\in k$ we have a homomorphism of complexes
${{\beta+t\cdot id_Z}\choose \psi}$. Let $$N_t:=ker({{\beta+t\cdot
id_Z}\choose \psi})$$ in the category of complexes. For any $t$
with $f_t:={{\beta+t\cdot id_Z}\choose \psi}$ being surjective we
have that $f_t$ is locally split. Here we call a homomorphism of
complexes $g$ locally split if $g$ is split in each degree, but
not necessarily split as a homomorphism of complexes. For all such
$t$ we see that $N_t$ is a complex of projective modules with the
same dimension array as $N$. We now consider $\pi_n(N_t)$,
$\pi_n(N)$, $\pi_n(M)$, $\pi_n(Z)$ and the induced mappings on the
truncated complexes. Of course, we still have
$ker(\pi_n(f_t))=\pi_n(N_t)$.

We shall prove that
$$t\mapsto \pi_n(N_t)\in \Scomp^{\underline{dim}{\pi_n(N)}}$$
is a rational morphism of varieties, imitating Christine Riedtmann's
proof in \cite{Riedtmann}.

There is an open neighborhood
$U$ of $0$ in $k$ so that $\pi_n(f_t)$ is surjective for all $t\in U$, using that
being surjective is an open condition and that $\pi_n(f_0)$ is
surjective.

Let $$(B,\partial^B)\stackrel{f}{\lra}
(A,\partial^A)\lra 0$$
be a surjective map of complexes of projective modules.
We want to compute
the kernel $(C,\partial^C)$ of this map. Since the structure of
$C$ as graded module is clear,
we may choose bases in $B$ so that we can identify
$B$ with $C\oplus
A$ as graded modules. Let $g=(g_C,g_A):C\lra B$ be the inclusion of the kernel.
We have $f={f_C\choose f_A}$
where $f_A$ is an isomorphism. Then
$g_C$ is an isomorphism
as well and we may assume that $g_C=id_C$. But then,
$$g_A=-f_Cf_A^{-1}$$
since $f_A$ is invertible. The differential on $ker(f)$ depending on $f$ is
$$
\partial_C =(id_C,-f_Cf_A^{-1})\cdot \partial_B\cdot{id_C\choose 0}\;.
$$
Thus we get a rational morphism of varieties
$Hom(B,A)\lra
\Scomp^{\underline{dim}{(C)}}$ defined on
the open neighborhood of
$f\in Hom(B,A)$ for which $f_A$
is an isomorphism.

We may now apply this construction to the map $f_0$ and by
composing with the map
$$
t\mapsto {\beta+t\cdot id_Z\choose \psi}
$$
we get the promised rational morphism of varieties.

Finally, for those $t$ for which $\pi_n(\beta+t\cdot id_Z)$ is an isomorphism,
that is for all but the finite number of eigenvalues of $-\beta$,
we get $\pi_n(N_t)\simeq \pi_n(M)$ and for $t=0$ we get
$\pi_n(N_0)\simeq \pi_n(N)$. Therefore
$$\pi_n(N)\in\overline{G\cdot\pi_n(M)}.$$
\phantom{x}\hfill\dickebox

\section{Geometric relation implies algebraic relation} \label{converse}

We shall prove in this section that under some conditions
the inverse implication of Theorem~\ref{adherence} is true as well.

Let $\underline d=(d_n,\dots,d_m)$ be a bounded dimension array.
We associate to the affine variety $\Scomp^{\underline d}(k)$
an affine $k$-scheme $\Scomp^{\underline d}(-)$. This $k$-scheme has
the following functorial description. Let $R$ be a commutative
$k$-algebra. Let $\Scomp^{\underline d}(R)$ denote the subset of
$$\prod_{i\in\Z}Hom_R(R^{\alpha(d_i)},R^{\alpha(d_{i-1})})$$
consisting of elements $(\partial_i)_{i\in \Z}$ with the
properties that $\partial _i$ is an $R\otimes _k A$-homomorphism
when viewed as a map from $\bigoplus_{j=1}^l (R\otimes_k
{P_j}^{{d^j_i}})$ to $\bigoplus_{j=1}^l (R\otimes_k
{P_j}^{{d^j_{i-1}}})$ and $\partial_i \partial_{i-1} = 0$. For a
$k$-algebra homomorphism $f:S\lra R$ there is naturally a
corresponding map $f^*:\Scomp^{\underline d}(S)\lra
\Scomp^{\underline d}(R)$ sending a tuple of matrices
$(\partial_i)$ to the tuple $(f(\partial_i))$. Similarly we may
associate to the affine algebraic group $G$ a smooth affine group
scheme $G(-)$ over $k$. The action of $G$ on $\Scomp^{\underline d}$
extends to an action of $G(-)$ on $\Scomp^{\underline d}(-)$.

We may verify Grunewald-O'Halloran's conditions which are necessary to apply
\cite[Theorem 1.2]{GH}.

\begin{Theorem}\label{geomtoalg}
Let $\underline d$ be a dimension array. Let $M,N\in \Scomp^{\underline d}$
be two complexes with bounded homology. If $M\leq_{top}N$
then $M\leq_{\Delta}N.$
\end{Theorem}

Proof. Let $n$ be an integer such that the homology of $M$ and $N$
vanishes in all degrees larger or equal to $n$. We will construct
a short exact sequence of complexes $$0\lra \pi_n(N) \lra \pi_n(M)
\oplus Z_{(n)} \lra Z_{(n)} \lra 0$$ where $Z_{(n)}$ is a complex of
projective $A$-modules. We are going to follow the steps of
Zwara's proof for the module case.

By Grunewald-O'Halloran's result \cite[Theorem 1.2]{GH}
there is a discrete valuation $k$-algebra $R$ with maximal ideal
${\mathfrak m}$ and residue field $k$ and with over $k$
finitely generated quotient field $K$ of transcendence degree $1$
and a complex $Y$ in $\Scomp^{\underline d_n}(R)$ so that
$k\otimes_RY=\pi_n(N)$ and as complexes of $K\otimes_RA$-modules,
$K\otimes_RY=g\cdot (K\otimes_k\pi_n(M))$ for a $g\in G(K)$. Since the
valuation on $R$ is discrete, $\mathfrak m$ is principal, generated by
an element $f$.

Since $\underline d_n$ is bounded there is
a non-zero element $z\in R$ so that $zg$ is a tuple of matrices
with entries in $R$. Using the
explicit definition of the action we get,
$$K\otimes_RY=g\cdot (K\otimes_k\pi_n(M))
=zg\cdot (K\otimes_k\pi_n(M)).$$ So, we may assume that $g$ is a tuple of matrices
with entries in $R$. Restricting
the multiplication with $g$ to $R\otimes_k\pi_n(M)$ gives a morphism of
complexes of $R\otimes_kA$-modules $\varphi:R\otimes_k\pi_n(M)\lra Y$.
Let $X$ denote the image of this morphism. Both $X$ and $Y$ are
complexes of free $R$-modules, with equal rank in all degrees,
therefore there exists some $s$ such that ${\mathfrak m}^sY\subseteq X$.

Now we take the point of view that the complexes $X$ and $Y$
are graded $R\otimes_kA$-modules with differentials. Fix a $k$-basis
$\mathcal{B}$ of $R$. As complexes of $A$-modules we have
$$X=\bigoplus_{b\in {\mathcal B}}X_b$$
where $X_b=\varphi(<b>\otimes _k\pi_n(M))\cong \pi_n(M)$ and where $<b>$ denotes
the $k$-subspace of $R$ generated by $b$.

For each $h$ we have a short exact sequence
of complexes $$0\lra Y/{\mathfrak m}Y\lra Y/{\mathfrak m}^{h+1}Y\lra
Y/{\mathfrak m}^{h}Y\lra 0$$ We will show that there exists an $h$ such that
$Y/{\mathfrak m}^{h+1}Y\simeq \pi_n(M)
\oplus \left(Y/{\mathfrak m}^{h}Y\right)$ as complexes of $A$-modules
where the mapping $\left(Y/{\mathfrak m}Y\right)\simeq
\left({\mathfrak m}^{h}Y/{\mathfrak m}^{h+1}Y\right)\lra
\left(Y/{\mathfrak m}^{h+1}Y\right)$ is induced by multiplication by $f^h$
and canonical inclusion. Let $V=\oplus_i V_i$ be a graded
vector space formed by taking vector space complements of
$X_{i}$ in $Y_{i}$ in each degree $i$.
Note that $V$ is a finite dimensional vector space since $Y$ is bounded and
${\mathfrak m}^{s}Y\subseteq X$.
Let $Z_0$ be the smallest $A$-subcomplex of $Y$ containing $V$.
Then $Z_0$ is a finite dimensional complex of $A$-modules,
since $Y$ is bounded. Now $Y = X + Z_0$ as complexes of $A$-modules.
Let $\mathcal{V}$ be a finite subset of $\mathcal{B}$ such that
$Z_0\cap \bigoplus_{b\in {\mathcal V}}X_b = Z_0\cap X$. Such a subset exists since
$Z_0$ is finite dimensional over $k$. Let $Z_1 = Z_0 + \bigoplus_{b\in {\mathcal V}}X_b$.
Then $Y = Z_1 \oplus \bigoplus_{b\not\in {\mathcal V}}X_b$. Since $\mathcal V$ is
finite there exists an integer
$t$ such that ${\mathfrak m}^{t+1}X \cap \bigoplus_{b\in {\mathcal V}}X_b=0$.
Thus there is a finite subset $\mathcal{W}$ of $\mathcal{B}$ such that
$${\mathfrak m}^{t+1}X \oplus \bigoplus_{b\in {\mathcal W}}X_b \oplus
\bigoplus_{b\in {\mathcal V}}X_b = X$$
Let $Z_2 = Z_1 + \bigoplus_{b\in \mathcal{W}}X_b$. Then $Y = {\mathfrak m}^{t+1}X \oplus Z_2$.

It follows that we have a chain of inclusions
$${\mathfrak m}^{s+t+2}Y\subseteq {\mathfrak m}^{t+2}X \subseteq
{\mathfrak m}^{t+1}X\subseteq Y$$ where the last two inclusions have
direct complements as complexes of $A$-modules. Thus
$$Y/{\mathfrak m}^{s+t+2}Y
\cong ({\mathfrak m}^{t+2}X/ {\mathfrak m}^{s+t+2}Y) \oplus \pi_n(M)
\oplus(Y/{\mathfrak m}^{t+1}X) \simeq \pi_n(M)
\oplus (Y/{\mathfrak m}^{s+t+1}Y)$$
where the last isomorphism follows since
$${\mathfrak m}^{t+2}X/{\mathfrak m}^{s+t+2}Y
\simeq {\mathfrak m}^{t+1}X/{\mathfrak m}^{s+t+1}Y\mbox{ and }
{\mathfrak m}^{t+1}X/{\mathfrak m}^{s+t+1}Y\oplus Y/{\mathfrak m}^{t+1}X
\simeq Y/{\mathfrak m}^{s+t+1}Y.$$
Now since
$Y/{\mathfrak m}Y\cong \pi_n(N)$ we get the promised short exact
sequence of complexes $0\lra \pi_n(N) \lra \pi_n(M)\oplus Z_{(n)} \lra
Z_{(n)} \lra 0$ by choosing $Z_{(n)} = Y/{\mathfrak m}^{s+t+1}Y$.

Now construct a complex $N'$ by splicing $\pi_n(N)$ with a projective
resolution $P_N$ of $H_n(\pi_n(N))$. Similarly we form a complex $Z$
by splicing a projective resolutions $P_Z$ of $H_n(Z_{(n)})$ with $Z_{(n)}$.
By the
horseshoe lemma there exists a short exact sequence $0\lra P_{N} \lra
P_{M\oplus Z} \lra P_{Z} \lra 0$ of projective resolutions where $P_{M\oplus Z}
\simeq P_N\oplus P_Z$ as graded modules and where $P_{M\oplus Z}$ is a
projective resolution of $H_n(\pi_n(M) \oplus Z_{(n)}) \cong H_n(\pi_n(M))
\oplus H_n(Z_{(n)})$. Moreover we have a short exact sequence of complexes
$$0\lra N' \lra M' \lra Z\lra 0$$ where the complex $M'$ is formed by
splicing $P_{M\oplus Z}$ with the complex $\pi_n(M)\oplus Z_{(n)}$. Now $N'$,
$M'$ are homotopy equivalent to $N$, $M\oplus Z$, respectively. Thus we
get a triangle $N\lra M\oplus Z \lra Z \lra N[1]$, which completes
the proof of the theorem. \hfill\dickebox

\section{Consequences for the geometry of complexes}
\label{consequences}

We continue with some consequences and observations on
$\Scomp^{\underline d}$ and the orders $\leq_{\Delta}$ and
$\leq_{top}$.

\begin{Example}
We consider the quiver $Q$ defined by $\bullet_1\lra\bullet_2$.
Then, up to isomorphism, there are $3$ indecomposable $kQ$-modules,
the indecomposable projective module $P_1 $ corresponding to the vertex
$1$ and the two simple modules $S_1$ and $P_2$. Moreover, in the
representation variety
$mod(kQ,(1,1))$ of $2$-dimensional $kQ$-modules with two different
composition factors, one has that the projective module with top $1$
degenerates to the direct sum of the two simple modules.
The projective indecomposable module with top $1$ can be considered
as being in $\Scomp^{\left({0\choose 0},{1\choose 0}\right)}$, where
$a\choose b$ indicates that in a certain degree the module is
$P_1^a\oplus P_2^b$. The semi-simple module $S_1\oplus S_2$ is in
$\Scomp^{\left({0\choose 1},{1\choose 1}\right)}$.
So, the modules are represented in different
varieties $\Scomp^{\underline d}$ and here it is not possible to consider
degenerations between them if one declares that
a complex $X$ degenerates to a complex $Y$ if $Y$ is in the closure of
the orbit of $X$.
Nevertheless, one may consider another non minimal projective
resolution of $P_1$ as
$$P_2\stackrel{(id_{P_2},0)}{\lra} P_2\oplus P_1.$$
This complex can be seen as being in
$\Scomp^{\left({0\choose 1},{1\choose 1}\right)}$,
and the minimal projective resolution of $S_1\oplus P_2$ is
$$P_2\stackrel{(0,\iota)}{\lra} P_2\oplus P_1.$$ for
$\iota$ being the embedding $P_2\lra P_1.$ Therefore, $P_1$ and
$S_1\oplus P_2$ can be both visualized in
$\Scomp^{\left({0\choose 1},{1\choose 1}\right)}$. Moreover it is easy
to see that $P_1\leq_{top} S_1\oplus P_2$.
This observation is one of the motivations not to ask for the complexes
to be minimal as is done in \cite{DrozdBekkert} and to allow zero
homotopic direct summands.
\end{Example}

Let $M$ and $N$ be $d$-dimensional $A$-modules. We write $M\leq N$
if $M$ degenerates to $N$ in $mod(A, d)$.

\begin{Prop}
Let $M,N\in mod(A, d)$ for some dimension $d$ and let
$P_M, P_N \in \Scomp^{\underline d}$ for some dimension array
${\underline d}$ be a projective resolution of $M$ and $N$, respectively.
Then, $M\leq N$ in $mod(A, d)$ if and only if
$P_M\leq_{top}P_N$  in $\Scomp^{\underline d}$.
\end{Prop}

Proof. If $M\leq N$, then by Zwara's theorem \cite{Zwaramain}
there is an exact sequence
$$0\lra N\lra Z\oplus M\lra Z\lra 0$$ for an $A$-module $Z$.
This implies a distinguished triangle
$$ P_N\lra P_Z\oplus P_M\lra P_Z\lra P_N[1]$$
in $K^-(A)$ where $P_Z$ is a projective resolution of $Z$. Hence,
$P_M\leq_{\Delta}P_N$ and so by Theorem~\ref{adherence} we have
$P_M\leq_{top}P_N$.

Conversely suppose $P_M\leq_{top}P_N$ and so by Theorem \ref{geomtoalg}
we have $P_M\leq_{\Delta}P_N$.
Then, there is a complex $Z$ and
a distinguished triangle $$P_N\lra Z\oplus P_M\lra Z\lra P_N[1]$$
Taking homology of this triangle gives a long exact sequence
$$
\lra H_{i+1}(Z)\lra H_i(P_N)\lra H_i(Z)\oplus H_i(P_M)\lra H_i(Z)\lra
H_{i-1}(P_N)\lra
$$
where $H_i(P_N)=H_i(P_M)=0$ for $i> 0$. For $i=0$ one gets
an exact sequence
$$0\lra H_{1}(Z)\lra H_{1}(Z)\lra N\lra H_0(Z)\oplus M\lra
H_0(Z)\lra 0\;.$$
This implies that  $$0 \lra N\lra H_0(Z)\oplus M\lra H_0(Z)\lra 0$$
is a short exact sequence and hence $M\leq N$ in
$mod(A,d)$, again by Zwara's theorem \cite{Zwaramain}.
This proves the statement. \hfill\dickebox

\bigskip


\begin{Lemma} Let $\underline d$ be any dimension array and
let $T$ be an element in $\Scomp^{\underline d}$ so that
$G\cdot T$ is open. Then, $T$ is a minimal element for $\leq_{\Delta}$
and for $\leq_{top}$.
\end{Lemma}

Proof.
If $G\cdot T$ is open, then $\Scomp^{\underline d}\setminus\{G\cdot T\}$
is closed and for any $X\not\simeq T$ one has that
$$
\overline{G\cdot X}\subseteq \Scomp^{\underline d}\setminus\{G\cdot T\}.
$$
Hence, $T$ is minimal with respect to $\leq_{top}$, and since
$\leq_{\Delta}$ implies $\leq_{top}$, the complex $T$ is minimal
also with respect to $\leq_{\Delta}$.
\phantom{x}\hfill\dickebox

\bigskip

Observe that we only used the topology of the space in the previous argument.
We shall see that for bounded $\underline d$ the orbits of $T$ with
$Hom_{D^b(A)}(T,T[1])=0$ are open.

\begin{Lemma} \label{openorbit}
Let $\underline d$ be a bounded dimension array and let $X$
be a complex in $\Scomp^{\underline d}$.
If $Hom_{D^b(A)}(X,X[1])=0$ then $G\cdot X$ is open in $\Scomp^{\underline d}$
\end{Lemma}

Proof.
First assume that ${\underline d}$ is a bounded dimension array.
From Theorem 7
in \cite{SaorinHuisgenZimmermann} we see that the orbit of $X$ in
$comp(A,\underline\alpha(\underline d))$ is open if
$Hom_{D^b(A)}(X,X[1])=0$.
The result now follows since
$\Scomp^{\underline d}$
is a subvariety of $comp(A,\underline\alpha(\underline d))$ and since
$G\cdot X = \Scomp^{\underline d}\cap
(Gl_{\underline\alpha(\underline d)} \cdot X)$. \hfill\dickebox

\bigskip

\begin{Lemma}\label{partorder}
The relation $\leq_{top}$ is a partial order on the set of
isomorphism classes of complexes with bounded homology with fixed dimension
array $\underline d$.
\end{Lemma}

Proof. If $N\in \overline{G\cdot M}$ and $M\in\overline{G\cdot L}$,
then clearly $N\in\overline{G\cdot L}$.
If $N\in\overline{G\cdot M}$ and $M\in\overline{G\cdot N}$, then
by the proof of Theorem~\ref{adherence} we get
$\pi_n(N)\in\overline{G\cdot \pi_n(M)}$ and
$\pi_n(M)\in\overline{G\cdot \pi_n(N)}$ for all $n\in\Z$. This implies
$\pi_n(N)\simeq\pi_n(M)$ for all $n\in\Z$

We show that whenever $X$ is a complex with bounded
homology in $\Scomp^{\underline d}$, then denoting by $m$
an integer so that the homology of $X$ is $0$ in all
degrees higher than $m$, then
$Y\in G\cdot X$ if and only if
$Y\in\pi_{\ell}^{-1}(G\cdot\pi_\ell(X))$
for all $\ell\geq m+1$. Indeed, assume that
$Y\in\pi_{\ell}^{-1}(G\cdot\pi_\ell(X))$ for
all $\ell\geq m+1$. Then we have an isomorphism of homology groups
$H_n(Y)\simeq H_n(X)$
for all $n$, which shows that $H_n(Y)=0$ for all $n\geq m+1$.
Then there is an isomorphism
$\pi_{m+1}(Y)\simeq \pi_{m+1}(X)$, which lifts to a homotopy equivalence
$Y\simeq X$ and
so $Y\in G\cdot X$. The reverse implication is trivial.

Hence, one has $N\simeq M$.\hfill\dickebox

\bigskip

\begin{Remark}
\cite[Theorem 7]{SaorinHuisgenZimmermann} cited in the proof
of Lemma~\ref{openorbit}
shows that the tangent space of the variety of complexes
$comp(A,\underline\alpha)$
at some point $X$ modulo the tangent space of the orbit of $X$ under the
group which is acting at $X$ is isomorphic to $Hom_{D^b(A)}(X,X[1])$.
A similar result can be proven for $\Scomp^{\underline d}$
and the action of our
smaller group. We also mention that Lemma \ref{openorbit}
has a converse in the
case where ${\underline d}$ is bounded. Namely, if $G\cdot X$ is open in
$\Scomp^{\underline d}$ then $Hom_{D^b(A)}(X,X[1])=0$.
This can again be seen
from \cite[Theorem~7]{SaorinHuisgenZimmermann}.
\end{Remark}

\begin{Cor}
Complexes with $Hom_{D^b(A)}(X,X[1])=0$ are minimal with respect to
both $\leq_{\Delta}$ and  $\leq_{top}$. In particular partial tilting
complexes are minimal with respect to both orders.
\end{Cor}

We also give a consequence which doesn't require an algebraically closed
field.

\begin{Cor} \label{unicity}
Let $A$ be an algebra over a field $K$. Then, up to homotopy equivalence
there is at most one two-term partial tilting complex
$$T=\dots\lra0\lra P_1\lra P_0\lra 0\lra\dots$$ with fixed homogeneous
components $P_0$ and $P_1$.
\end{Cor}

Proof. Since two-term complexes of projective modules are entirely
determined by their homology, and since for any field extension $L$
of $K$ one has $H(L\otimes_KX)\simeq L\otimes_KH(X)$ for any complex
$X$, we may assume that $K$ is algebraically closed. Let
$\alpha_i=\underline{dim}(P_i)$ for $i\in\{0,1\}$ and
$\underline\alpha:=(\alpha_1,\alpha_0)$. The variety $\Scomp^{\underline
d(\underline\alpha)}$ is an affine space, and therefore irreducible
as algebraic variety. Moreover, since $T$ is a partial tilting complex,
the orbit $G\cdot T$ is open in $\Scomp^{\underline d(\underline\alpha)}$.
Therefore, $G\cdot T$ is dense. Let $S$ be another partial tilting complex in
$\Scomp^{\underline d(\underline\alpha)}$. Also $G\cdot S$ is open and
dense, and therefore $S\leq_{top}T$ as well as $T\leq_{top}S$. Hence,
$S\simeq T$ by Lemma~\ref{partorder}.  \hfill\dickebox

\bigskip

\begin{Example}
Corollary \ref{unicity} does not hold for general dimension arrays.
Let $A$ be given by the quiver
$$\unitlength1cm
\begin{picture}(1.8,.85)
\put(0,.35){$\bullet_1$}
\put(.2,.6){\vector(1,0){1.3}}
\put(1.4,.25){\vector(-1,0){1.3}}
\put(1.5,.35){$\bullet_2$}
\put(.6,.65){\scriptsize $\alpha$}
\put(.6,.0){\scriptsize $\beta$}
\end{picture}$$
with relations $\alpha\beta\alpha=\beta\alpha\beta=0$.
For this algebra take the indecomposable complex (unique up to
isomorphism so that $P_2$ is in degree $0$)
$$T_1:=\dots\lra 0\lra P_1\lra P_1\lra P_2\lra 0\lra\dots$$
and $T_2:=P_1\lra P_2$. Then, $T:=T_1\oplus T_2$ is a tilting complex.
Let $\underline d$ be the dimension array of $T$.
The complex $S$
$$P_1\stackrel{(0,id)}{\lra}P_1\oplus
 P_1\stackrel{\left(\begin{array}{cc}
\alpha&0\\0&0\end{array}\right)}{\lra} P_2\oplus P_2$$
is homotopy equivalent to the tilting complex $S'$
$$P_1\stackrel{(\alpha, 0)}{\lra} P_2\oplus P_2.$$
Here both tilting complexes $T$ and $S$ have the same dimension array,
but are not isomorphic, and therefore belonging to different irreducible
components of $\Scomp^{\underline d}$. Using \cite{Lakshmibai}
and a slightly more detailed examination of $\Scomp^{\underline d}$ one
observes that $\Scomp^{\underline d}$ has exactly two irreducible
components.

The complex $T_1\oplus P_1[2]$ is a tilting complex as well and denote by
$\underline e$ the dimension array of $T_1\oplus P_1[2]$. A short
examination yields that $\Scomp^{\underline e}$ has two irreducible
components, one $C_3$ of dimension $3$ and another
component $C_4$ of dimension $4$. The orbit of $T_1\oplus P_1[2]$ is
open in $C_3$, whereas the complexes
corresponding to the points in $C_4$ are not partial tilting complexes.
Observe however, in $C_4$ there is an open orbit of a complex $U\simeq
P_1[2]\oplus P_2$ with $Hom_{D^b(A)}(U,U[1])=0\neq Hom_{D^b(A)}(U,U[2])$.
\end{Example}

\begin{Rem}
Observe that a tilting complex $T$ over $A$ is the image $F(B)$ of an
equivalence $F:D^b(B)\lra D^b(A)$ of triangulated categories.
By the Rickard's and Keller's
main theorem \cite{Ri3,Keller2} there is a so-called twosided
tilting complex $X$ of $A\otimes_KB^{op}$-modules which are projective
on the left and on the right, so that $X\otimes_B^{\mathbb L} -$
is an equivalence.
For any dimension array $\underline d$ let $X\otimes \underline d$ be
the dimension array which is obtained by tensoring a complex with
dimension array $\underline d$ by $X$, and taking the total complex of
the resulting bi-complex. Then, by definition
$X\otimes_B-$ induces a  morphism of varieties
$$\Scomp(X):\Scomp^{\underline d}_B\lra \Scomp^{X\otimes\underline d}_A\;.$$
It should be an interesting question to study the image of this morphism
inside $\Scomp^{X\otimes\underline d}_A$. Note that studying varieties
using functors is already far from trivial in the module case (see Bongartz
\cite{BongartzMorita} and Zwara \cite{Zwarafunctor}).
\end{Rem}

There is another consequence of these statements. Indeed, define
for any two complexes $X$ and $Y$ in $K^{-,b}(A)$
$$X\leq_{Hom}Y:\Leftrightarrow \forall{U\in D^b(A)}:\;
dim_k(Hom_{D^b(A)}(U,X))\leq dim_k(Hom_{D^b(A)}(U,Y)) $$

\begin{Lemma}
Let $X$ and $Y$ be two complexes in $\Scomp^{\underline d}$
for bounded dimension array $\underline d$. Then,
$X\leq_{top}Y\Rightarrow X\leq_{Hom}Y$.
\end{Lemma}

Proof. Define for any two complexes $X$ and $Y$ with appropriate
bounded dimension array $\underline d$ and $\underline e$ the mapping
$$\varphi_{X,Y}:\prod_{i\in\Z}Hom_A(X_i,Y_{i+1})\lra Hom_{C^b(A)}(X,Y)$$
by $\varphi_{X,Y}(f):=\partial_Xf+f\partial_Y.$ It is clear that
this image is exactly the set of $0$-homotopic homomorphisms.
Hence, we have that
$$dim_k(Hom_{D^b(A)}(X,Y))=dim_k(Hom_{C^b(A)}(X,Y))-
dim_k(im(\varphi_{X,Y}))\;.$$
We use the argument from \cite[\S 3, Theorem 2, special case]{CBlecture}
to show that $$\{U\}\times \Scomp^{\underline d} \lra \N$$
given by $(U,X)\mapsto dim_k(Hom_{D^b(A)}(U,X))$ is upper
semi-continuous. Then, setting $n=dim_k(Hom_{D^b(A)}(U,X)),$ one gets
$\{Z| dim_k(Hom_{D^b(A)}(U,Z))\geq n\}$ is closed, and if
$Y\in\overline{G\cdot X}$, then $Y\in \{Z| dim_k(Hom_{D^b(A)}(U,Z))\geq
n\}.$ Hence, $$dim_k(Hom_{D^b(A)}(U,Y))\geq dim_k(Hom_{D^b(A)}(U,X)).$$
This proves the statement.\hfill\dickebox

\bigskip

\paragraph{\bf Acknowledgment:}
We would like to thank the referee for helpful remarks.

\end{document}